\documentclass[12pt,leqno]{article}
\usepackage{amssymb}
\newcommand{\dd}{{\rm \kern 3pt I\kern-9pt d}}

\newcommand{\Abar}{{\backslash\kern-8pt A}}

\title{\large ON THE COARSEST TOPOLOGY PRESERVING CONTINUITY}
\author{\sc Nicolas Bouleau\footnote{ENPC, ParisTech, 28 rue des Saints P\`eres, 75007 Paris; 
e-mail : {\tt bouleau@enpc.fr}}}
\date{juin 2006}

\font\helv Helvetica at 14pt
\newcommand{\CC}{\mbox{{\helv \i}{\rm\kern-5pt C}}}

\begin{document}
\maketitle

\noindent{\bf R\'esum\'e.} Nous \'etudions sur un espace de fonctions une topologie, appel\'ee topologie collante, qui a la propri\'et\'e d'\^etre la moins fine parmi celles qui pr\'eservent la continuit\'e. Consid\'er\'ee dans les cadres convenables, cette topologie conserve le caract\`ere bor\'elien, l'int\'egrabilit\'e locale, le caract\`ere c\`adl\`ag et d'autres propri\'et\'es. Elle est moins fine que la convergence uniforme locale et elle s'accommode du ph\'enom\`ene de bosse glissante ainsi que nous le montrons sur des exemples. Nous \'etablissons des crit\`eres de relative compacit\'e pour cette topologie et nous envisageons quelques extensions.\\

\noindent{\bf Abstract.} We study a topology on a space of functions, called sticking topology, with the property to be the weakest among the topologies preserving continuity. In suitable frameworks, this topology preserves borelianity, local integrability, right continuity and other properties. It is coarser than the locally uniform convergence and it allows the presence of gliding humps as we show on examples. We prove relative compactness criteria for this topology and we consider some extensions.\\

\noindent{\Large\textsf{I. Introduction.}}\\

There exists a topology between the point-wise and the locally uniform convergences, which shares with the second one the property to preserve continuity. We call it the topology of sticky convergence or sticking topology. It is the coarsest topology preserving continuity in the sense that its restriction to the set of continuous functions coincides with the point-wise convergence.

We study some features of this topology and show that it preserves several properties generally damaged by point-wise convergence (borelianity, local $\mu$-integrability, right continuity, c\`adl\`ag or l\`adl\`ad functions, etc.). Let us say right now that to the natural question whether this topology may yield efficient tools to provide weak convergence criteria for laws of stochastic processes, up to now the answer is  disappointing. The sticking topology seems to be too closed to the point-wise convergence to bring not yet available significant information. For instance, a set of continuous functions is relatively compact in the set of all functions for this topology if and only if it is relatively compact for the point-wise convergence and any cluster point in continuous.

In order to be as concrete as possible and to give quickly the main ideas, we consider first the case of $\mathbb{R}^{\mathbb{R}_+}$ the space of real functions on $[0,+\infty[$ (Part II). The sticky convergence is defined by its neighbourhoods. As the point-wise convergence it is not metrizable. We prove its relations with point-wise convergence and its properties with respect to continuity.

The next part is devoted in the same setting to the properties preserved by the sticky convergence. Several qualitative properties are preserved, some quantitative properties yield continuous or lower semi-continuous functionals.

Part IV presents some examples. The peculiar feature of sticky convergence is to allow situations usually considered as pathological and to be excepted. A sequence $f_n$ may converge for the sticking topology to a constant although the number of up-crossings $N^{a,b}(f_n)$ of the interval $[a,b]$ becomes unbounded. This is due to the possibility of gliding humps which prevent equicontinuity but allow the continuity of the limit. Our aim here is just to give some non obvious examples of the gliding humps phenomenon which is for a long time a central concept in analysis. It appears especially in connection with Schur's lemma on the equivalence between weak and strong convergence for sequences in $\ell^1$ (cf. [12] p 122), with the Hellinger-Toeplitz theorem saying that infinite matrices sending $\ell^2$ in $\ell^2$ are norm continuous and with the Vitali-Hahn-Sachs theorem. Gliding humps techniques have generated specific arguments gathered in a genuine theory (cf. [11]) that we dont tackle here.

The fifth part restores all their generallity to the ideas by showing that the sticky convergence may be defined  on the space of applications from a topological space to an other topological space. It comes from a uniform structure when the second space is uniform. We take this opportunity (prop. 9) to correct an omission in [1].

Some ways of constructing relatively compact sets for the sticking topology are given in part VI. We introduce in particular the notion of a set ``equicontinuous after the humps" which possesses similar features as an equicontinuous set but allows the presence of humps. On the space of c\`adl\`ag functions the Skorohod topology shares with the sticking topology the property to be coarser than the locally uniform convergence but it is not comparable with point-wise convergence. The same happens for the pseudo-paths topology studied by Meyer and Zheng [9] and the S-topology introduced by Jakubowski [7]. On the other hand, these topologies possess sufficient (Meyer-Zheng) or necessary and sufficient (Jakubowski) criteria of relative compactness well adapted to the stochastic calculus on semi-martingales (see e.g. in [6] comments to this question).

The last part proposes some extensions of the peculiar idea of the sticky convergence: proximity on a neighbourhood after possible humps. On $\ell^\infty$ the sticky convergence at infinity defines a metrizable topological vector space. By isomorphism thanks to a Schauder basis or directly, some topologies may be defined on Banach spaces allowing humps.

With respect to [1] where appeared the initial idea, the present article attempts to go further to make this kind of arguments more popular.\\

\noindent{\Large\textsf{II. The sticking topology.}}\\

In the whole article topological spaces are supposed to be Hausdorff (in particular uniformizable spaces are completely regular). We consider first the case of the space $X= \mathbb{R}^{\mathbb{R}_+}$ of functions from $[0,\infty[$ into $\mathbb{R}$.

An elementary neighbourhood of $f\in X$ is defined as a set of the form
\begin{equation}
E^\varepsilon_{t_1,\ldots,t_k}(f)=\{g : \exists \eta>0,\;\forall t\in\cup_{i=1}^k]t_i-\eta,t_i+\eta[\quad|f(t)-g(t)|<\varepsilon\}
\end{equation}
An open set is a set containing an elementary neighbourhood of each of its points. This defines a topology called topology of sticky convergence or sticking topology. Comparing neighbourhoods  gives immediately :\\

\noindent{\bf Proposition 1.} {\it The sticking topology is Hausdorff, finer than the point-wise convergence and coarser than the locally uniform convergence.}\\

A net $(f_\alpha)$, $\alpha$ in a directed set $(D,\succ)$ (cf. Choquet [4]), converges sticky to $f$ if 
$$
\forall t\in\mathbb{R}_+,\forall\varepsilon>0, \exists\alpha\in D, \forall\beta\succ\alpha, \exists V_t {\mbox{ neighbourhood of }} t, \forall s\in V_t\;|f_\beta(s)-f(s)|<\varepsilon.
$$

\noindent{\bf Proposition 2.} {\it a) The set $\mathcal{C}([0,\infty[,\mathbb{R})$ of continuous functions is closed in $X$ equipped with the sticking topology.

\noindent b) On $\mathcal{C}([0,\infty[,\mathbb{R})$ the sticking and point-wise topologies are equivalent.}\\

\noindent{\bf Proof.} a) If $h$ belongs to the sticky closure of $\mathcal{C}([0,\infty[,\mathbb{R})$, in every $E^\varepsilon_{t_1}(h)$ there are continuous functions. Let $g$ be some of them, there exists 
$\eta_1>0$ such that $|t-t_1|<\eta_1\Rightarrow|h(t)-g(t)|<\varepsilon.$ 
Now $g$ being continuous there exists $\eta_2>0$ such that $|t-t_1|<\eta_2\Rightarrow|g(t)-g(t_1)|<\varepsilon.$ 
So that for $ |t-t_1| < \eta_1\wedge \eta_2$ we have 
$|h(t)-h(t_1)|<2\varepsilon$, proving $h$ is continuous at $t$.

b) Let $f_\alpha\in\mathcal{C}([0,\infty[,\mathbb{R})$ be a net point-wise converging to the continuous function $f$. Let us consider a sticky neighbourhood $E^\varepsilon_{t_1,\ldots t_k}(f)$ of $f$. By the point-wise convergence there is an $\alpha$ such that $\max_{i=1,\ldots,k}|f_\alpha(t_i)-f(t_i)|<\varepsilon/3$. By the fact that both $f$ and $f_\alpha$ are continuous, there is $\eta>0$ such that $\forall t\in\cup_i]t_i-\eta,t_i+\eta[$ holds $|f_\alpha(t)-f(t)|<\varepsilon$.  This proves $f_\alpha\in E^\varepsilon_{t_1,\ldots t_k}(f)$.\hfill$\Box$

The following remarks up to the end of this part, might be done for filters, sequences are just used for the sake of simplicity. There exists for the sticking topology a criterion of convergence which, as the Cauchy criterion, does not involve the limit. Let us denote $A^o$ the interior of the set $A$.\\

\noindent{\bf Proposition 3.} {\it A sequence $(f_n)$ converges sticky iff}
\begin{equation}
\forall t, \forall\varepsilon>0, \exists N,\forall n>N, \exists V_t, \forall s\in V_t, \exists M,\forall m>M,\;|f_m(s)-f_n(s)|<\varepsilon.
\end{equation}
{\it where $V_t$ is a neighbourhood of $t$. This property may also be written}
\begin{equation}
\forall\varepsilon>0,\; \liminf_n(\liminf_m\{s: |f_m(s)-f_n(s)|<\varepsilon\})^o=\mathbb{R}_+.
\end{equation}

The proof is quite similar to that of proposition 2. Let us remark that if we drop the interior operator $(.)^o$ in (3) we obtain a weaker condition
\begin{equation}
\forall\varepsilon>0,\; \liminf_n\liminf_m\{s: |f_m(s)-f_n(s)|<\varepsilon\}=\mathbb{R}_+
\end{equation}
which says that $\forall s\in \mathbb{R}_+$ $(f_n(s))$ is a Cauchy sequence, i.e. point-wise convergence. We can observe also that the locally uniform convergence may be written
\begin{equation}
\cup_N(\cap_{n\geq N}\cup_M\cap_{m\geq M}\{u : |f_m(u)-f_n(u)|<\varepsilon\})^o=\mathbb{R}_+
\end{equation}
which is indeed  stronger than (3) since $(\cap_nA_n)^o\subset\cap_nA_n^o$.\\

\noindent{\Large\textsf{III. Properties preserved by the sticky convergence.}}\\

A property $P$ is said to be preserved if the set of $f$ in $X$ satisfying $P$ is closed, i.e.
$$
\left(\forall t,\forall\varepsilon>0, \exists V_t, \exists g\in P, (|f-g|<\varepsilon\mbox{ on }V_t)\right)\Rightarrow f\in P.
$$
We saw already that continuity is preserved.\\

\noindent{\bf Proposition 4.}{\it Let $\tau=(t_k)$ be a sequence of real numbers tending to $t$. The functionals $S_\tau(f)=\limsup_k f(t_k)$ and $I_\tau(f)=\liminf_k f(t_k)$ with values in $\overline{\mathbb{R}}$ are  continuous for the sticking topology.}\\

\noindent{\bf Proof.} A net $(f_\alpha)$ which converges sticky to $f$ satisfies
$$\forall\varepsilon>0, \exists\alpha, \forall\beta\succ\alpha, \exists V_t :|f_\beta-f|<\varepsilon\mbox{ on }V_t.
$$
It follows that
$$\limsup_k f(t_k)-\varepsilon\leq \limsup_k f_\beta(t_k)\leq\limsup_k f(t_k)+\varepsilon
$$ hence $\lim_\alpha S_\tau(f_\alpha)$ exists and is equal to $S_\tau(f)$ and similarly for $I_\tau$.\hfill$\Box$\\

\noindent{\bf Proposition 5.}{\it The following properties are preserved :

a) lower and upper semi-continuity,

b) right continuity, to be c\`ag, c\`ad, l\`ag, l\`ad, c\`adl\`ag, etc.

c) local boundedness,

d) borelianity, $\mathcal{F}$-measurability ($\mathcal{F}$ $\sigma$-field on $\mathbb{R}_+$ containing the Borel sets), local $\mu$-integrability ($\mu$ measure on $\mathbb{R}_+$), approximate continuity at $a\in\mathbb{R}_+$,

e) continuity (resp. c\`ad, etc.) on a dense $G_\delta$ set, continuity (resp. c\`ad, etc.) outside a $\mu$-negligible set.}\\

\noindent{\bf Demonstration.} The points a) b) c) are similar to the argument for continuity. The properties of point d) are in general related to sequences, they are valid here because the starting space $\mathbb{R}_+$ is locally compact. Let us consider for instance the $\mathcal{F}$-measurability. Let $f$ be such that
$$\forall t, \forall\varepsilon>0, \exists V_t, \exists g\; \mathcal{F}\mbox{-measurable  and }(|f-g|<\varepsilon \mbox{ on }V_t).
$$
Let $K$ be a compact set in $\mathbb{R}_+$, $K$ being covered by finitely many $V_t$'s, there exists $g_{K,\varepsilon}$ $\mathcal{F}$-measurable such that $|f-g_{K,\varepsilon}|<\varepsilon$ on $K$, hence $f$ is $\mathcal{F}$-measurable.

Let us recall that $h$  approximately continuous at $t$ means that there exists $A\subset\mathbb{R}_+$ Lebesgue negligible such that $h(a)=\lim_{x\rightarrow a, x\in\hspace{-0.15cm}/A} h(x)$. If $h_\alpha$ is a net converging sticky to $h$, by the above argument there exists a countable family $(h_\beta)\subset(h_\alpha)$ such that $h$ belongs to the sticky closure $\overline{\{h_\beta\}}^s$ of $\{h_\beta\}$. Then, if $A_\beta$ are the corresponding Lebesgue negligible sets, taking $A=\cup_\beta A_\beta$, the property follows from the preservation of continuity.

Similarly if $f\in\overline{A}^s$ there is a countable subset of $A$, say $A_0$, such that $f\in\overline{A_0}^s$. Then the property concerning $G_\delta$ sets comes from the fact that $\mathbb{R}_+$ is a Baire space (cf. Choquet [4]).\hfill$\Box$\\

\noindent{\bf Remarks.} {\bf 1.} Right continuity is preserved the topology of $\mathbb{R}_+$ being the usual one. This is of course also true when $\mathbb{R}_+$ is equipped with the right topology whose open sets are generated by intervals $[a,b[$, the new sticking topology is coarser than the preceding one.

{\bf 2.} Let $(\Omega, \mathcal{A}, \mathbb{P})$ be a probability space with an increasing family of $\sigma$-fields $\mathcal{F}_t$, and let $X_n(t)$ be a sequence of real stochastic processes. If $\mathbb{P}$-almost surely $X_n$ converges sticky to $X$, (or even if there are modifications of a subsequence of the $X_n$'s which converge sticky to $X$) then, if the $X_n$'s are respectively adapted, c\`adl\`ag, progressively measurable, then $X$ has almost surely the same property.

{\bf 3.} Let us consider the number of up-crossings of the interval $[a,b]$ on $[0,1]$.
$$N^{a,b}(f)=\max\{k : 0\leq t_1\leq\cdots\leq t_{2k}\leq1\quad f(t_{2i-1})<a<b<f(t_{2i})\;\forall i\}$$ the property $(N^{a,b}(f)<\infty\;\;\forall a<b)$ which means that $f$ is l\`adl\`ag is preserved. The functional $N^{a,b}(f)$ is lower semi-continuous for the point-wise convergence and therefore also for the sticking convergence. Now when the $f_n$ are continuous and $f_n\rightarrow f$ sticky, the number $N^{a,b}$ doesn't remain bounded in general.\\

For comparison with the locally uniform convergence we have\\

\noindent{\bf Proposition 6.} {\it a) Let $(f_n)$ be a sequence of continuous functions. If there are continuous functions $g_n$ which converge locally uniformly to $g_\infty$ and if 
$$\limsup_n\{f_n\neq g_n\}=\emptyset$$ then the sequence $(f_n)$ converges sticky to $g_\infty$.

b) If the $f_n$'s are no more supposed to be continuous, with the same hypotheses on the $g_n$'s, if 
$$\limsup_n\overline{\{f_n\neq g_n\}}=\emptyset$$ then $f_n$ converges sticky to $g_\infty$.}\\

\noindent{\bf Proof.} a) The hypothesis means that for every $t$, $f_n(t)=g_n(t)$ for sufficiently large $n$. So $\lim_n f_n(t)=g_\infty(t)$ and $g_\infty$ is continuous, the result follows from prop 2.

b) Taking the opposite sentence, the hypothesis means that for every $t$, $t$ is in all the sets $\{f_n=g_n\}^o$ except possibly finitely many, hence 
$$\exists N, \forall n\geq N, \exists V_t \;(f_n(s)=g_n(s)\mbox{ on }V_t).$$ Now, if $V_t^\prime$ is a relatively compact neighbourhood of $t$, 
$$\exists M, \forall m\geq N, \forall s\in V_t^\prime \; \;|g_n(s)-g_\infty(s)|<\varepsilon$$
so that on $V_t\cap V_t^\prime$ we have $|f_n-g_\infty|<\epsilon$, and $f_n\rightarrow g_\infty$ sticky.\hfill$\Box$\\

\noindent{\Large\textsf{IV. Examples.}}\\

Sticky convergence without local uniformity will be essentially related to humps phenomena.\\

\noindent{\bf 1.} Let $\xi$ be a continuous function on $\mathbb{R}_+$ vanishing at zero and at infinity. The sequence $(f_n)$ defined by $f_n(t)=\xi(nt)$ converges sticky to zero.

Let us consider a signal $\eta(t)$ perturbed by infinitely many humps $$\eta_n(t)=\eta(t)+\sum_k\alpha_k\xi_k(n(t-t_k))$$ where the functions $\xi_k$ satisfy the same hypotheses as $\xi$, and $t_k$ are points in $\mathbb{R}_+$. 

Supposing $|\xi_k|\leq 1$ and $\sum_k|\alpha_k|<+\infty$, as easily seen using the dominated convergence in the sum, we have $\eta_n\rightarrow\eta$ sticky as $n\rightarrow\infty$. The convergence is not locally uniform already when there is only one hump.\\

\noindent{\bf 2.} Let $\xi\in\mathcal{S}(\mathbb{R})$ the set of infinitely differentiable rapidly decreasing functions, and suppose $\xi(0)=0$ and $\int_\mathbb{R}\xi(s)ds=0$.

If $\psi_N(t)=\sum_{n=-N}^N\xi(nt)$ then $\psi_N\rightarrow\psi_\infty$ sticky and $\psi_\infty$ is continuous at zero.\\

Indeed, if $t\neq 0$ the convergence is uniform on a neighbourhood of $t$. In $t=0$, let us put $\hat{\xi}(u)=\int e^{iux}\xi(x)dx$. The Poisson summation formula gives for $s\neq 0$
$$\sum_n\xi(sn)=\frac{1}{s}\sum_m\hat{\xi}(\frac{2\pi m}{s}).$$
Letting $s\rightarrow 0$, since $\hat{\xi}(0)=0$, and $\hat{\xi}$ being rapidly dedreasing $\sum_n\xi(ns)\rightarrow 0$. The hypotheses on $\xi$ may be of course weakened.\hfill$\Box$\\

\noindent{\bf 3.} Let $\eta$ be a function such that $\int_0^1\eta(s)ds=0$ and vanishing outside the unit interval. We put 
$$\eta_n(t)=\alpha_n\;n\,\eta(nt)$$
with $|\alpha_n|\rightarrow +\infty$ and we consider the convolution operators $R_n$ defined by
$$R_n(f)=f\star \eta_n\quad\quad f\in\mathcal{C}(\mathbb{T}).$$ On the Baire space $\mathcal{C}(\mathbb{T})$ the norm of $R_n$ is $|\alpha_n|\int_0^1|\eta(s)|ds$, thus by the Banach-Steinhaus theorem there are $f\in\mathcal{C}(\mathbb{T})$ such that $R_n(f)$ does not remain bounded.

Let us show on an example that  among these functions, there are in general some such that $R_n(f)\rightarrow 0$ point-wise.

Let us take $\eta(t)=2(1_{[0,\frac{1}{2}[}-1_{[\frac{1}{2},1[})$ and define $\zeta_k^{t_0}(t)= 0$ outside $]t_0-\frac{1}{k}, t_0+\frac{1}{k}[$ and $\zeta_k^{t_0}(t_0)=k$, $\zeta_k^{t_0}$ being affine on $[t_0-\frac{1}{k}, t_0]$ and on $[t_0,t_0+\frac{1}{k}]$. Then as easily verified

\noindent{\bf Lemma.} {\it For $n\geq k\geq4,$ $ \|\zeta^{t_0}_k\star\eta_n\|_\infty=\frac{k^2}{2n}|\alpha_n|$. This maximum is reached at $t_0$ and the support of $\zeta_k^{t_0}\star\eta_n$ is contained in $[t_0-\frac{1}{k},t_0+\frac{1}{k}+\frac{1}{n}].$}

We construct functions $f_n$ by induction putting $f_0=\beta_0\zeta^{\frac{1}{2}}_4$ and $f_i=\beta_i\zeta^{t_i}_{k_i}$ in such a way that the points $t_i\downarrow 0$ and the supports of the $f_i$ do not overlap and be separated by open intervals. Then we put $f=\sum_i f_i$. We have $f\in\mathcal{C}(\mathbb{T})$ as soon as $\beta_ik_i\rightarrow0$, we have also $f\star\eta_n(0)=0,$ and $f\star\eta_n=\sum_if_i\star\eta_n$ goes to zero at every $t\neq0$ as soon as $\frac{|\alpha_n|}{n}\rightarrow0$.

We can proceed in such a way that $(f\star\eta_n)(t_i)=\beta_i\alpha_n\frac{k^2_i}{2n}$ and $\beta_ik_i\rightarrow0$, $\frac{\alpha_n}{n}\rightarrow0$ allows to have nevertheless $\|f\star\eta_n\|_\infty\rightarrow+\infty.$\hfill$\Box$

Similarly, for Fourier series, if we put $S_n(f)=\sum_{-n}^n\hat{f}(k)e^{2i\pi kt}=\int_\mathbb{T}D_n(t-s)f(s)ds$ the Dirichlet kernels $D_n$ are not bounded in $L^1(\mathbb{T})$ and there are continuous functions whose Fourier series $S_n(f)$ converges to $f$ point-wise, $\|S_n(f)\|_\infty$ being unbounded.\\

\noindent{\bf Remark 4.} When $f_n\rightarrow f$ sticky, $f_n$ continuous as in the preceding examples, if $t_k$ is a sequence converging to $t_\infty$, the point $P= f(t_\infty)$ has the property to be a {\it double cluster point} to the double sequence $\sigma_{ij}=f_i(t_j)$, in the sense that every neighbourhood $V$ of $P$ meets infinitely many columns along infinitely many points and infinitely many rows along infinitely many points.

From the sticky convergence we have indeed
$$\exists N, \forall n>N, \exists K,\forall k>K,\;f_n(t_k)\in V.
$$ what shows that after $N$ the row $\sigma_{nj}$ meets $V$ for $j$ sufficiently large. Now $f$ is continuous and $f_n\rightarrow f$ point-wise, so that
$$\exists K,\forall k>K, \exists N, \forall n>N,\; f_n(t_k)\in V$$
and this shows that after $K$ the column $\sigma_{ik}$ meets $V$ fot $i$ sufficiently large.

We shall come back to this remark in Part VI.\\

\noindent{\Large\textsf{V. General setting.}}\\

Let $E$ and $F$ be topological spaces and $\mathcal{F}(E,F)$ the set of applications from $E$ to  $F$. A basis of neighbourhoods of $f\in\mathcal{F}(E,F)$ for the sticking topology is given by 
$\mathcal{V}^{G_1,\ldots,G_k}_{x_1,\ldots,x_k}$, $x_i\in E$, $G_i$ neighbourhoods of $f(x_i)$ in $F$, defined by
$$
\mathcal{V}^{G_1,\ldots,G_k}_{x_1,\ldots,x_k}=\{g : \exists V_{x_1},\ldots,V_{x_k}\mbox{ neighbhds of }x_1,\ldots,x_k\mbox{ resp. s. t. }(g(x)\in G_i \;\forall x\in V_{x_i})\}
$$
So defined the sticking topology is finer than the point-wise convergence. If the topology of $F$ is induced by a uniform structure with entourages $W$, the sticking topology is induced by the uniform structure whose entourages are 
$$\mathcal{W}^W_{x_1,\ldots,x_k}=\{(f,g):\exists V_{x_1},\ldots,V_{x_k}\mbox{ s. t. }((f(x),g(x))\in W\;\forall x\in\cup_i V_{x_i})\}.
$$
It is coarser than the locally uniform convergence.

The space $\mathcal{C}(E,F)$ of continuous functions is sticky closed and the restriction of sticky convergence is equivalent on $\mathcal{C}(E,F)$ with that of point-wise convergence. It follows that when $F$ is a complete uniform space, $\mathcal{F}_s(E,F)$, i.e. $\mathcal{F}(E,F)$ with the sticking topology, is never complete except in obvious cases where $\mathcal{C}(E,F)$ is point-wise closed.\\

\noindent{\bf Proposition 7.} {\it Suppose $F$ be a complete uniform space. Let $(D,\succ)$ be a directed set and $(f_\alpha)_{\alpha\in D}$ a net. In order that $(f_\alpha)$ be sticky convergent it is necessary and sufficient that
$$
\begin{array}{c}
\forall W\mbox{ entourage of } F, \forall x\in E, \exists\alpha\in D, \forall\beta\succ\alpha\\
\exists V_x\mbox{ neighbourhood of }x\mbox{ s. t. }\exists\gamma\in D,\forall\delta\succ\gamma\\
(f_\beta,f_\gamma)\in W
\end{array}
$$}
We refer to [1] prop 5. for the proof. When $F$ is only sequentially complete (e.g. a reflexive Banach space equipped with the weak topology, the space of radon measures on a locally compact space, the space of distributions) a necessary and sufficient condition holds for sequences.\\

\noindent{\bf Proposition 8.} {\it Let $F$ be sequentially complete. A sequence $(f_n)$ is sticky convergent iff
\begin{equation}
\begin{array}{l}
\forall W\mbox{ entourage of } F, \forall x\in E, \exists N,\forall n\geq N,\\ \exists V_x\mbox{ neighbourhood of }x\mbox{ s. t. }\exists M,\forall m\geq M, \;(f_m,f_n)\in W \mbox{ on } V_x.
\end{array}
\end{equation}}\\

As a consequence of this criterion, we have\\

\noindent{\bf Proposition 9.} {\it Let $(f_n)$ be a sequence of continuous functions from $E$ into a Banach space $(B,\|.\|)$.

a) If the series $\sum_k\|f_n(x)\|$ converges to a continuous function, the series $\sum_n f_n(x)$ converges in $B$ to a continuous function.

b) If $\|\sum_{k=0}^nf_k(x)\|\leq A(x)$, $A(x)$ being locally bounded, if $\varepsilon_n(x)$ is a sequence of continuous real functions converging point-wise to zero and such that $\sum_n|\varepsilon_{n+1}(x)-\varepsilon_n(x)|$ converges to a continuous function, then $\sum_n\varepsilon_n(x)f_n(x)$ converges to a continuous function.

c) If $f_n$ converges to a continuous function and if  $(g_n)$ is a sequence of continuous functions satisfying
$$\|g_p(x)-g_q(x)\|\leq K(x)\|f_p(x)-f_q(x)\|\quad\forall p,q>N(x)\in\mathbb{N}
$$where $K$ is locally bounded, then $g_n$ converges to a continuous function.}\\

\noindent{\bf Proof.} The point a) comes from the criterion and the inequality $\|\sum_p^qf_n(x)\|\leq \sum_p^q\|f_n(x)\|$. The point b) is proved as the Abel rule for real series. For c) let $V_x^o$ be a neighbourhood of $x$ where $K(y)\leq K_0$, and let $\varepsilon^\prime=\frac{\varepsilon}{K_0}$. Then 
$$\exists N,\forall p>N, \exists V_x, \forall y\in V_x, \exists M,\forall q>M\quad \|f_p(y)-f_q(y)\|<\varepsilon^\prime.$$ On $V_x\cap V_x^o$ we have therefore $\|g_p(y)-g_q(y)\|<\varepsilon$ and the criterion applies to the $g_n$'s.

\hfill$\Box$\\

\noindent{\Large\textsf{VI. Relative compactness.}}\\

Let $\mathcal{F}_s(E,F)$ [resp. $\mathcal{F}_p(E,F)$] denote the space $\mathcal{F}(E,F)$ equipped with the sticking [resp. point-wise] topology.\\

\noindent{\bf 1. Compactness, continuous functions and sequences.}\\

\noindent{\bf Proposition 10.}Ê {\it  Let be $H\subset\mathcal{C}(E,F)$. $H$ is relatively compact in $\mathcal{F}_s(E,F)$ iff

a) $\forall x\; H(x)$ is relatively compact in $F$,

b) $\overline{H}^p\subset\mathcal{C}(E,F)$,

\noindent where $\overline{H}^p$ denotes the closure of $H$ in $\mathcal{F}_p(E,F)$.}\\

\noindent{\bf Demonstration.} 1) The condition is sufficient because the restrictions of $\mathcal{F}_s(E,F)$ and  $\mathcal{F}_p(E,F)$ coincide on $\mathcal{C}(E,F)$.

2) Let $H$ be relatively compact in $\mathcal{F}_s(E,F)$. The map $f\mapsto f(x)$ is continuous from $\mathcal{F}_s(E,F)$ to $F$, hence a) is fulfilled. Let $\overline{H}^s$ be the closure of $H$ in $\mathcal{F}_s(E,F)$, $\overline{H}^s\subset\overline{H}^p$ and since the identity map from $\mathcal{F}_s(E,F)$ into $\mathcal{F}_p(E,F)$ is continuous, $\overline{H}^s$ is compact in $\mathcal{F}_p(E,F)$, hence closed, then $H\subset\overline{H}^s\subset\overline{H}^p$ gives $\overline{H}^s=\overline{H}^p$, so that $\overline{H}^p\subset\mathcal{C}(E,F)$.

\hfill$\Box$

\noindent{\bf Remark 5.} Let us note that the subsets $H$ of $\mathcal{C}(E,F)$ satisfying conditions a) and b)  of prop 10 are exacly the sets which in the topological space ($\mathcal{C}(E,F)$, point-wise) are relatively compact.\\

\noindent{\bf Proposition 11.} {\it Let $E$ be a countable union of compact sets ($E$ is a $K_\sigma$), and suppose $F$ metrisable. Then if $f_0$ is a cluster point of $A$ in $\mathcal{F}_s(E,F)$ there exists a countable subset $A_0$ of $A$ such that $f_0\in\overline{A_0}^s$.}\\

\noindent{\bf Demonstration.} Let $d$ be a distance generating the structure of $F$. Let $K$ be a compact subset of $E$. For $m,n\in \mathbb{N}$, let us show that there exists a finite part $A^K_{m,n}$ of $A$ such that
\begin{equation}
\forall x_1,\ldots,x_m\in K, \exists V_{x_1},\ldots, V_{x_m},\; \exists f\in A^K_{m,n}\;:\quad d(f_0,f)\leq \frac{1}{n}\mbox{ on } \cup_iV_{x_i}.
\end{equation}
By $f_0\in\overline{A}^s$ we have
\begin{equation}
\begin{array}{l}
\forall(t_1,\ldots,t_m)\in K^m,\exists \mbox{ open neighbourhoods } V_{t_1},\ldots,V_{t_m}\\
\exists f\in A\mbox{ s. t. }\;d(f_0,f)\leq \frac{1}{n}\mbox{ on } \cup_iV_{t_i}.
\end{array}
\end{equation}
The open sets $V_{t_1}\times\cdots\times V_{t_m}$ for $(t_1,\ldots, t_m)\in K^m$ cover the compact set $K^m$. There exists a finite sub-covering denoted 
$$G_i=V_{s_{i,1}}\times\cdots\times V_{s_{i,m}}\quad i=1,\ldots,N$$
thus $\forall(t_1,\ldots, t_m)\in K^m$ there is $G_i$ with $(t_1,\ldots, t_m)\in G_i$ and then there is $f\in A$ such that $d(f_0,f)\leq \frac{1}{n}$ on $\cup_iV_{s_{i,j}}$. The $N$ functions involved constitute a finite set $A^K_{m,n}$ satisfying (7).

Now since $E$ is a $K_\sigma$ the conclusion follows.\hfill$\Box$\\

\noindent{\bf Remark 6.} If $E$ is compact, since Dirac masses are particular cases of measures, on $\mathcal{C}(E,\mathbb{R})$ the weak topology $\sigma(\mathcal{C}(E,\mathbb{R}),\mathcal{M}(E))$ is finer than the point-wise topology. Actually if $H$ satisfies conditions a) and b) of prop. 10 and is {\it uniformly bounded}, then $H$ is weakly relatively compact ([5] thm 5). On such a set the weak, sticking and point-wise topologies coincide. 

The uniformly bounded subsets of the first Baire class behave also well with respect to measures (cf. [10]).\\
\newpage
\noindent{\bf 2. Relatively compact subsets of $\mathcal{F}_s(E,F)$.}\\

Let us come back to the case where $E$ is a general topological space. 

If $F$ is a topological vector space, it is easily seen that if $A$ and $B$ are relatively compact in $\mathcal{F}_s(E,F)$ then so is $A+B$.
In particular if $A\subset\mathcal{F}(E,F)$ is equicontinuous and if $A(x)$ is relatively compact in $F$ for any $x\in E$, then if $(f_n)$ is a sticky convergent sequence, the set 
$$B=\{f+f_n : f\in A, n\in \mathbb{N}\}$$
is sticky relatively compact. Therefore, if the $f_n$ are continuous [resp. c\`ad if $E=\mathbb{R}$, etc.] any sequence (or filter) in $B$ possesses a cluster point in $\mathcal{F}_s(E,F)$ which is continuous [resp. c\`ad, etc.], hence any point-wise convergence sequence (or filter) in $B$ converges to a continuous [resp. c\`ad, etc.] limit.

If we focuse on continuous functions, a criterion of relative compactness in $\mathcal{F}_s(E,F)$ reduces to conditions on a set $H$ in order that any point-wise limit of functions in $H$ be continuous.

For this the {\it Eberlein-Grothendieck criterion} is sufficient in a quite general framework and also most often necessary (cf. [5]).\\

\noindent{\bf Definition 1.} {\it A set $H\subset\mathcal{C}(E,F)$ satisfies the Eberlein-Grothendieck criterion if there exists $E_1$ dense in $E$ such that for any sequence $(f_n)$ in $H$ and any sequence $(x_i)$ in $E_1$, there exists $y\in F$ double cluster point of the double sequence $f_n(x_i)$, in the sense that any neighbourhood of $y$ meets infinitely many columns each along infinitely many points and the same for rows.}\\

\noindent{\bf Prposition 12.}(Grothendieck) {\it If $H\subset\mathcal{C}(E,F)$ satisfies the Eberlein-Grothendieck criterion, then any point-wise limit (sequence or filter) of elements of $H$ is continuous.}\\

(Therefore if $H(x)$ is relatively compact for any $x\in E$, $H$ is relatively compact in $\mathcal{F}_s(E,F)$.)

\noindent{\bf Demonstration.} It is remarkable that by the properties of completely regular spaces (Hausdorff uniformizable spaces) the problem may be reduced to the case where $F=\mathbb{R}$, we refer to Grothendieck [5] for this reduction.

If $f\in\overline{H}^p$ is not continuous at $x_0\in E$, there is $\varepsilon>0$ such that for any neighbourhood $V$ of $x_0$ there is an $x\in V\cap E_1$ such that $|f(x_0)-f(x)|\geq\varepsilon$. Then by induction we can construct a sequence of functions $(f_i)$ in $H$ and a sequence $(x_j)$ in $E_1$ such that
$$\begin{array}{ll}
a)\quad& | f_n(x_i)-f(x_i)|\leq \frac{1}{n}\quad\forall i=0,\ldots, n-1\\
b)\quad& |f_i(x_n)-f_i(x_0)|\leq \frac{1}{n}\quad\forall i=0,\ldots,n\\
c)\quad& |f(x_n)-f(x_0)|\geq \varepsilon
\end{array}
$$
now this is easily seen to be in contradiction with the existence of a  double cluster point to the double sequence $f_i(x_j)$.\hfill$\Box$\\

We can express the criterion in a different way. Let us come back to the case $E=\mathbb{R}_+$, $F=\mathbb{R}$, for the sake of simplicity.\\

\noindent{\bf Definition 2.} {\it A double real sequence $\sigma_{ij}$ will be said ``flat on the edges" is there is an application $\kappa$ from $\mathbb{N}$ to $\mathbb{N}$ such that the following limit exists
$$\lim_{\begin{array}{c}
i,j\rightarrow\infty\\
|i-j|\geq \kappa(i\wedge j)
\end{array}}
\sigma_{ij}
$$}

As easily seen the limit of a double sequence flat on the edges is a double cluster point for this sequence, so we have:\\

\noindent{\bf Corollary.} {\it Let be $H\subset\mathcal{C}(E,F)$, if for any sequence $(f_n)$ in $H$ and any sequence $(x_m)$ in $E$, the double sequence $f_i(x_j)$ is flat on the edges, then any point-wise limit of $H$ is continuous.}\\

The Eberlein-Grothendieck criterion uses a property to be fulfilled by any sequence $(f_n)$ in $H$ and says roughly that it cannot point-wise converge to a discontinuous function. A more quantitative sufficient criterion easier to verify is the following:\\

For $C$ countable in $\mathcal{F}(E,F)$ we define for $s,t\in\mathbb{R}_+$
$$\rho_C(s,t)=\inf_{
\begin{array}{c}
\mathcal{J}\mbox{ finite}\\
\mathcal{J}\subset C
\end{array}}
\sup_{f\in C\backslash\mathcal{J}}
|f(s)-f(t)|
$$
with the convention $\sup\emptyset=0$, and for any $A\subset\mathcal{F}(E,F)$ we put
$$
e_A(s,t)=\sup_{
\begin{array}{c}
C\mbox{ countable}\\
C\subset A
\end{array}}
\rho(s,t)
$$

\noindent{\bf Definition 3.} {\it $H\subset\mathcal{F}(E,F)$ will be said ``equicontinuous after the humps" if $\forall t\in E$ ther exists $B\subset H$ equicontinuous at $t$ such that
$$\lim_{s\rightarrow t} e_{H\backslash B}(s,t)=0.$$}

\noindent{\bf Proposition 13.} {\it Let be $H\subset\mathcal{C}(\mathbb{R}_+,\mathbb{R})$ equicontinuous after the humps, then if $H(t)$ is relatively compact for each $t$, then $H$ is relatively compact in $\mathcal{F}_s(E,F)$.}\\

\noindent{\bf Proof.} Let us first remark that $e_A(s,t)$ is increasing in $A$ and that if $H$ is equicontinuous after the humps, so is any subset of $H$.

By prop. 10 and 11, it is sufficient to show that any sequence in $H$ possesses a continuous cluster point. Let $(f_n)$ a sequence in $H$ and $f$ a point-wise cluster point of $(f_n)$. The subset $S=\{f_n\}$ is equicontinuous after the humps. Let $t\in E$, there is  $B$ equicontinuous at $t$ such that $\lim_{s\rightarrow t}e_{S\backslash B}(s,t)=0$.

a) If $f$ is a cluster point of $B$ then $f$ is continuous at $t$.

b) If $f$ is a cluster point of $S\backslash B=\{f_{n_k}\}$ by point-wise convergence
$$\begin{array}{rl}
|f(s)-f(t)|&\leq \limsup_k|f_{n_k}(s)-f_{n_k}(t)|\\
&=\rho_C(s,t)\quad\mbox{ for }C=\{f_{n_k}\}\\
&\leq e_{S\backslash B}(s,t)\rightarrow0
\end{array}
$$
and $f$ is also continuous at $t$ in this case.\hfill$\Box$\\

If we restrict the setting to the Skorohod space $D$ of c\`adl\`ag functions which is closed in $\mathcal{F}_s(E,F)$, a net converges sticky iff it converges sticky on the rational numbers. It follows that the vector space $D$ equipped with the sticking topology is metrizable (but not Polish).

As a consequence if for any $\varepsilon>0$, probability measures $\mathbb{P}_n$ are carried up to $\varepsilon$ by a compact set $K_\varepsilon\subset D$ and if the finite marginal laws converge to $\mathbb{P}$, then $\mathbb{P}_n$ converge weakly to $\mathbb{P}$ for the sticking topology, hence by a classical argument
$$\int \Phi\,d\mathbb{P}_n\rightarrow\int\Phi\,d\mathbb{P}$$
for any bounded functional $\Phi$ from $D$ into $\mathbb{R}$ whose set of dicontinuity points are $\mathbb{P}$-negligible.

But, the price to be paid for the fact that the sticking topology is the coarsest preserving continuity is that there are very few such functionals, essentially \break $\Phi(f)=\lim_{t_k\rightarrow t}f(t_k)$ or finite continuous combination of them.\\

\noindent{\Large\textsf{VII. Extensions.}}\\

Let us add some comments around the preceding properties. The simplest case where the sticking topology may be considered is on the set $\ell^\infty$ of bounded real sequences. For $u\in\ell^\infty$ let us put $$\|u\|_s=\sum_k\frac{1}{2^k}|u(k)|+\limsup_k|u(k)|$$ and define $\ell_s$ as the space of bounded sequences equipped with the norm $\|.\|_s$ (sticking convergence at infinity). As already seen $\ell_s$ is not complete but a sequence converges in $\ell_s$ iff
$\|\limsup_m|u_m-u_n|\|_s\rightarrow 0$ as $n\rightarrow\infty$ (cf. prop 3). The classical spaces $\ell^p, 0<p\leq\infty, c_0,c$ are closed in $\ell_s$. The dual space of $\ell_s$ contains $\ell^1$ and other functionals obtained by the Hahn-Banach theorem.
Let us recall that a Banach limit is a linear form $\phi$ on $\ell^\infty$ such that

i) $\phi(u)=\phi(\theta(u))$ where $\theta(u)(k)=u(k+1)$

ii) $\liminf_k u(k)\leq\phi(u)\leq\limsup_k u(k)$

\noindent by ii) such a $\phi$ is continuous on $\ell_s$.

The same happens for medial limits defined by Mokobodzki thanks to the continuum hypothesis (cf [8]).\\

Starting with $\ell_s$ it is possible to transport the topology to Banach spaces using a Schauder basis. For instance through the historical Schauder basis in $\mathcal{C}([0,1])$ we obtain a topology on $\mathcal{C}([0,1])$ allowing humps (whose supports are getting smaller and smaller) in which several usual spaces (H\"{o}lder spaces, etc.) are closed.

Now, directly without using $\ell_s$, if we consider the space of bounded functions in an open domain $\Omega\subset\mathbb{R}^d$ with boundary $\partial\Omega$, a metrizable uniform structure may be defined by saying that a sequence $f_n$ converges if the convergence is locally uniform inside $\Omega$ and for a decreasing sequence $U_p$ of neighbourhoods of $\partial\Omega$ 
$$\limsup_p\|(f_n-f)|_{U_p}\|_\infty\rightarrow0\quad\mbox{ as }{n\rightarrow\infty}.$$
In this setting the space of continuous functions up to the boundary is closed and humps going toward the boundary are allowed.

\begin{list}{}
{\setlength{\itemsep}{0cm}\setlength{\leftmargin}{0.5cm}\setlength{\parsep}{0cm}\setlength{\listparindent}{-0.5cm}}
  \item\begin{center}
{\small REFERENCES}
\end{center}\vspace{0.4cm}
[1] {\sc Bouleau N.} ``Une structure uniforme sur un espace $\mathcal{F}(E,F)$." {\it Cahiers Topologie G\'eom. Diff.} Vol XI, 2, (1969), 207-214.

[2] {\sc Bourbaki N.} {\it Topologie G\'en\'erale}, Chap. 2 Structures uniformes, Hermann 1965.

[3] {\sc Bourbaki N.} {\it Espaces Vectoriels Topologiques}, Hermann 1967.

[4] {\sc Choquet G.} {\it Lectures on Analysis}, t I,II,III, Benjamin 1969.

[5] {\sc Grothendieck A.} ``Crit\`eres de compacit\'e dans les espaces g\'en\'eraux" {\it Amer. J. Math.} Vol 74, (1952), 168-186.

[6] {\sc Jacod J.; Protter Ph.} ``Asymptotic error distribution for the Euler method for stochastic differential equations" {\it Annals Prob.} Vol 26, 1? (1998), 267-307.

[9] {\sc Meyer P. A.; Zheng W. A.} ``Tightness criteria for laws of semi-martingales" {\it Ann. Inst. Henri Poincar\'e} Vol 20, 4, (1984), 353-372.

[10] {\sc Rosenthal H.} ``Point-wise compact subsets of the first Baire class" {\it Amer. J. Math. } Vol 99, 2, (1977), 362-378.

[11] {\sc Swartz Ch.} {\it Infinites matrices and the gliding hump}, World Scientific 1997.

[12] {\sc Yosida K.} {\it Functional Analysis}, Springer 1974.

\end{list}

\end{document}